\documentclass[11pt,reqno]{amsart}
\usepackage{amssymb,amsmath}

\newcommand{\be}{\begin{eqnarray}}
\newcommand{\ee}{\end{eqnarray}}

\newcommand{\eps}{{\mbox{$\epsilon$}}}

\newcommand{\R}{{\mathbb R}}

\newcommand{\C}{{\mathbb C}}

\newcommand{\LL}{\mathcal L}

\newcommand{\G}{{\mathcal G}}

\newtheorem{theorem}{Theorem}\newtheorem{lemma}[theorem]{Lemma}

\theoremstyle{definition}
\theoremstyle{remark}\numberwithin{equation}{section}\input epsf.sty
\begin{document}\thispagestyle{empty}

\title{Buffon needle lands in $\epsilon$-neighborhood of a $1$-Dimensional Sierpinski Gasket with probability at most $|\log\epsilon |^{-c}$}
\author{Matt Bond}\address{Matt Bond, Department of Mathematics, Michigan State University,
{\tt bondmatt@msu.edu}}
\author{Alexander Volberg}\address{Alexander Volberg, Department of Mathematics, Michigan State University,
{\tt volberg@math.msu.edu}}

%\thanks{Research of the authors was supported in part by NSF grants  DMS-0501067 (Nazarov and Volberg) and DMS-0605166 (Peres)}

\maketitle

\begin{abstract}

In recent years, relatively sharp quantitative results in the spirit of the Besicovitch projection theorem have been obtained for self-similar sets by studying the $L^p$ norms  of the ``projection multiplicity" functions, $f_\theta$, where $f_\theta(x)$ is the number of connected components of the partial fractal set that orthogonally project in the $\theta$ direction to cover $x$. In \cite{NPV}, it was shown that $n$-th partial 4-corner Cantor set with self-similar scaling factor $1/4$ decays in Favard length at least as fast as $\frac{C}{n^p}$, for $p<1/6$. In \cite{BV}, this same estimate was proved for the $1$-dimensional Sierpinski gasket for some $p>0$. A few observations were needed to adapt the approach of \cite{NPV} to the gasket: we sketch them here. We also formulate a result about all self-similar sets of dimension $1$.

\end{abstract}

%\newcommand{\K}{{\mathcal K}}
%\newcommand{\G}{{\mathcal G}}
%\newcommand{\LL}{{\mathcal L}}
%newcommand{\eps}{{\epsilon}}

\section{Definitions and result}
\label{result}
Let $E\subset\C$, and let $\text{proj}_\theta$ denote orthogonal projection onto the line having angle $\theta$ with the real axis. The \textbf{average projected length} or \textbf{Favard length} of $E$, $\text{Fav}(E)$, is given by
$$\text{Fav}(E)=\frac1{\pi}\int_0^{\pi}|\text{proj}_\theta(E)|d\,\theta.$$
For bounded sets, Favard length is also called \textbf{Buffon needle probability}, since up to a normalization constant, it is the likelihood that a long needle dropped with independent, uniformly distributed orientation and distance from the origin will intersect the set somewhere.

$B(z_0,r):=\{z\in\C:|z-z_0|<r\}$. For $\alpha\in\lbrace -1, 0, 1\rbrace^n$ let $$z_\alpha:=\sum_{k=1}^n{(\frac{1}{3})^ke^{i\pi[\frac{1}{2}+\frac{2}{3}\alpha_k]}},\,\,\,\G_n:=\bigcup_{\alpha\in\lbrace -1, 0, 1\rbrace^n} B(z_\alpha,3^{-n}).$$
This set is our approximation of a partial Sierpinski gasket; it is strictly larger. We may still speak of the approximating discs as ``Sierpinski triangles."

The main result:
\begin{theorem}\label{main}
$\text{Fav}(\G_n)\leq\frac{C}{n^{c}}\,,\, c>0$.
\end{theorem}
Set $\G_n$ is $3^{-n}$ approximation to Besicovitch irregular set (see \cite{falc1} for definition) called Sierpinski gasket. Recently one detects a considerable interest in estimating the Favard length of such $\epsilon$-neighborhoods of Besicovitch irregular sets, see \cite{PS}, \cite{T}, \cite{NPV}, \cite{LZ}. In \cite{PS} a random model of such Cantor set is considered and estimate $\asymp \frac1n$ is proved. But for non-random self-similar  sets the estimates of \cite{PS} are more in terms  of $\frac1{\log\cdots\log n}$ (number of logarithms depending on $n$) and more suitable for general class of ``quantitatively Besicovitch irregular  sets" treated in \cite{T}.

Let $f_{n,\theta}:=\frac1{2}\nu_n *3^n\chi_{[-3^{-n},3^{-n}]},$ where
$$\nu_n :=*_{k=1}^n\widetilde{\nu}_k\text{ and }\widetilde{\nu}_k:=\frac{1}{3}[\delta_{3^{-k}cos(\pi/2 -\theta)} +\delta_{3^{-k}cos(-\pi/6 -\theta)} +\delta_{3^{-k}cos(7\pi/6 -\theta)}].$$

For $K>0$, let $A_K:= A_{K,n,\theta}:=\{x:f_{n,\theta}\geq K\}$. Let $\LL_{\theta,n}:=\text{proj}_\theta(\G_n)=A_{1,n,\theta}$.
For our result, some maximal versions of these are needed:
$$f_{N,\theta}^*:=max_{n\leq N}f_{n,\theta},\,\,\, A_K^*:= A_{K,n,\theta}^*:=\{x:f_{n,\theta}^*\geq K\}.$$

Also, let $E:=E_N:=\{\theta:|A_K^*|\leq K^{-3}\}$ for $K=N^{\epsilon_0}, \,\epsilon_0$.

Later, we will jump to the Fourier side, where the function $$\varphi_\theta(x):=\frac1{3}[e^{-i\cos(\pi/2-\theta)}+e^{-i\cos(-\pi/6-\theta)}+e^{-i\cos(7\pi/6-\theta)}]$$
plays the central role: $\widehat{\nu_n}(x)=\prod_{k=1}^n\varphi_\theta(3^{-k}x)$.

\section{General philosophy}
\label{phil}

Fix $\theta$. If the mass of $f_{n,\theta}$ is concentrated on a small set, then $||f_{n,\theta}||_p$ should be large for $p>1$ - and vice versa. $1=\int f\leq ||f_{n,\theta}||_p||\chi_{\LL_{\theta,n}}||_q$, so $m(\LL_{\theta,n})\geq ||f||_p^{-q}$, a decent estimate. The other basic estimate is not so sharp: 
\begin{equation}\label{naive}m(\LL_{\theta,N})\leq 1-(K-1)m(A_{K,N,\theta})\end{equation}
However, a combinatorial self-similarity argument of \cite{NPV} and revisited in \cite{BV} shows that for the Favard length problem, it bootstraps well under further iterations of the similarity maps:

\begin{theorem}
\label{}
If $\theta\notin E_N$, then $|\LL_{\theta,NK^3}|\leq\frac{C}{K}$.
\end{theorem}

Note that the maximal version $f_N^*$ is used here. A stack of $K$ triangles at stage $n$ generally accounts for more stacking per step the smaller $n$ is. For fixed $x\in A_{K,N,\theta}^*$, the above theorem considers the smallest $n$ such that $x\in A_{K,n,\theta}$, and uses self-similarity and the Hardy-Littlewood theorem to prove its claim by successively refining an estimate in the spirit of \eqref{naive}.
Of course, now Theorem \ref{main} follows from the following:
\begin{theorem}\label{mofE}
Let $\epsilon_0<1/11$. Then for $N>>1$, $|E|<N^{-\epsilon_0}$.
\end{theorem}
It turns out that $L^2$ theory on the Fourier side is of great use here. It is proved in \cite{NPV}, \cite{BV}:
\begin{theorem}\label{L2}
For all $\theta\in E_N$ and for all $n\leq N$, $||f_{n,\theta}||_{L^2}^2\leq CK$.
\end{theorem}
One can then take small sample integrals on the Fourier side and look for lower bounds as well.
Let $K=N^{\epsilon_0}$, and let $m=2\epsilon_0 \log_3N$. Theorem \ref{L2} easily implies the existence of $\tilde{E}\subset E$ such that $|\tilde{E}|>|E/2|$ and number $n$, $N/4<n<N/2$, such that for all $\theta\in\tilde{E}$,
$$\int_{3^{n-m}}^{3^n}{\prod_{k=0}^n{|\varphi_\theta(3^{-k}x)|^2}dx}\leq \frac{2CKm}{N}\leq 2\eps_0N^{\eps_0-1}\log N.$$
Number $n$ does not depend on $\theta$; $n$ can be chosen to satisfy the estimate in the average over $\theta\in E$, and then one chooses $\tilde{E}$. Let $I:=[3^{n-m},3^n].$

Now the main result amounts to this (with absolute constant $A$ large enough):
\begin{theorem}\label{lower}
$$\theta\in\tilde{E}: \int_I{\prod_{k=0}^n{|\varphi_\theta(3^{-k}x)|^2}dx}\geq c3^{m-2\cdot Am}=cN^{-2\epsilon_0(2A-1)}.$$
\end{theorem}
The result: $2\eps_0 \log N\geq N^{1-\eps_0(4A-1)}$, i.e., $N\leq N^*$.
Now we sketch the proof of Theorem \ref{lower}. We split up the product into two parts: high and low-frequency: $P_{1,\theta}(z)=\prod_{k=0}^{n-m-1}\varphi_{\theta}(3^{-k}z)$, $ P_{2,\theta}(z)=\prod_{k=n-m}^n \varphi_{\theta}(3^{-k}z)$.

\begin{theorem}\label{P1}
For all $\theta\in E$, $\int_I|P_{1,\theta}|^2\,dx\geq C\,3^m\,.$
\end{theorem}
 Low frequency terms do not have as much regularity, so we must control the damage caused by the \textbf{set of small values}, $SSV(\theta):=\{x\in I:|P_2(x)|\leq 3^{-\ell}\}$, $\ell=\alpha\,m$ with sufficiently large constant $\alpha$. In the next result we claim the existence of $\mathcal{E}\subset\tilde{E}$, $|\mathcal{E}|>|\tilde{E}/2|$ with the following property:

\begin{theorem}\label{P1SSV}
$$\int_{\tilde{E}}\int_{SSV(\theta)}|P_{1,\theta}(x)|^2dx\,d\theta\leq 3^{2m-\ell/2}\Rightarrow \forall \theta\in \mathcal{E}\,\,\int_{SSV(\theta)}|P_{1,\theta}(x)|^2dx\,d\theta\leq c\,K\,3^{2m-\ell/2}\,.$$
\end{theorem}
 Then Theorems \ref{P1} and \ref{P1SSV} give Theorem \ref{lower}.

\section{Locating zeros of $P_2$}
\label{zP2}
We can consider $\Phi(x, y) = 1+ e^{ix} +e^{iy}$. The key observations are
$$
|\Phi(x,y)|^2 \ge a(|4\cos^2\,x -1|^2 +|4\cos^2 \,y-1|^2)\,,\,\,\,\frac{\sin 3x}{\sin x} = 4\cos^2 \,x -1\,.
$$
Changing variable we can replace $3\varphi_{\theta} (x)$ by $\phi_t(x)=\Phi (x, tx)$. Consider $P_{2,t}(x):=\prod_{k=n-m}^{n} \frac13\phi_t(3^{-k} x)$, $P_{1,t}(x):=\prod_{k=0}^{n-m} \frac13\phi_t(3^{-k} x)$. We need $SSV(t):=\{x\in I: |P_{2,t}(x)|\le 3^{-\ell}\}$. One can easily imagine it if one considers $\Omega:=\{(x,y)\in [0,2\pi]^2: |\mathcal{P}(x,y)|:=|\prod_{k=0}^m \Phi( 3^k x, 3^k y) |\le 3^{m-\ell}\}$. Moreover, (using that if $x\in SSV(t)$ then $3^{-n}x\ge 3^{-m}$, and using $xdxdt=dxdy$) we change variable in the next integral:
$$
\int_{\tilde{E}}\int_{SSV(t)} |P_{1,t}(x)|^2\,dxdt =3^{-2n+2m}\cdot 3^n\int_{\tilde{E}}\int_{3^{-n}SSV(t)}|\prod_{k=m}^n \Phi( 3^k x, 3^k tx)|^2\,dxdt \le 
$$
$$
3^{-n+3m}\int_{\Omega} |\prod_{k=m}^n \Phi( 3^k x, 3^k y)|^2\,dxdy\,.
$$
Now notice that by our key observations $\Omega\subset \{(x,y)\in [0,2\pi]^2: |\sin 3^{m+1} x|^2 +|\sin 3^{m+1}y|^2 \le a^{-m} 3^{2m-2\ell}\le  3^{-\ell}\}\,.$ The latter set $\mathcal{Q}$ is the union of $4\cdot 3^{2m+2}$ squares $Q$ of size $3^{-m-\ell/2}\times 3^{-m-\ell/2}$. Fix such a $Q$ and estimate
$$
\int_{Q} |\prod_{k=m}^n \Phi( 3^k x, 3^k y)|^2\,dxdy\le 3^{\ell} \int_{Q} |\prod_{k=m+\ell/2}^n \Phi( 3^k x, 3^k y)|^2\,dxdy \le 
$$
$$
3^{\ell} \cdot (3^{-m-\ell/2})^2\int_{[0,2\pi]^2}|\prod_{k=0}^{n-m-\ell/2} \Phi( 3^k x, 3^k y)|^2\,dxdy \le
3^{\ell} \cdot (3^{-m-\ell/2})^2\cdot 3^{n-m-\ell/2}= 3^{-2m} \cdot 3^{n-m-\ell/2}\,.
$$
Therefore, taking into account the number of squares $Q$ in $\mathcal{Q}$ and the previous estimates we get

$$
\int_E\int_{SSV(t)} |P_{1,t}(x)|^2\,dxdt\le 3^{2m-\ell/2}\,.
$$
Theorem \ref{P1SSV} is proved.

To prove Theorem \ref{P1} 
we need the following simple lemma.
\begin{lemma}
\label{CETSQ}
Let $C$ be large enough. Let $j=1,2,...k$, $c_j\in\C$, $|c_j|=1$, and $\alpha_j\in\R$. Let $A:=\lbrace \alpha_j\rbrace_{j=1}^k$. Suppose
$$\int_{\R} (\sum_{\alpha\in A} \chi_{[\alpha-1, \alpha+1]}(x))^2\,dx \le S\,.\text{ Then }\int_0^1 |\sum_{\alpha\in A} c_{\alpha} e^{i\alpha y}|^2\, dy \le C\,S\,.$$
\end{lemma}
Some key facts useful for its proof:
$$
\int_0^1|\sum_{\alpha\in A} c_{\alpha} e^{i\alpha\,y}\, dy|^2 \le e\int_0^{\infty}|\sum_{\alpha\in A} c_{\alpha} e^{i(\alpha+i)\,y}\, dy|^2=e\int_{\R} \bigg|\sum_{\alpha\in A}\frac{c_{\alpha}}{\alpha +i -x}\bigg|^2\,dx\, ,
$$
and the fact that $H^2(\C_+)$ is orthogonal to $\overline{H^2(\C_+)}$, so one can pass to the Poisson kernel.

\section{The general case}
\label{general}

Let us have $k$ closed disjoint discs of radii $1/k$ located in the unit disc. We build $k^n$ small discs of radii $k^{-n}$ by iterating $k$ linear maps from small discs onto the unit disc. Call the resulting union $S_k(n)$. We would like to show that exactly as in the case of $k=3$ considered above  and in a very special case of $k=4$ considered in \cite{NPV} $\text{Fav}(S_k(n))\le C\, n^{-c}, \, c>0$. However, presently we can prove only a weaker result.

\begin{theorem}
\label{gent}
$$
\text{Fav}(S_k(n))\le C\, e^{-c\,(\log n)^{1/2}}, \, c>0\,.
$$
\end{theorem}

\markboth{}{\sc \hfill \underline{References}\qquad}

\end{document}